\documentclass[a4paper,11pt,oneside,reqno]{amsart}
\usepackage{amsmath,amssymb,amsthm,mathtools}
\usepackage[pagewise]{lineno} 
\usepackage{graphicx}
\usepackage{color}
\mathtoolsset{showonlyrefs}
\usepackage{bm}

\usepackage{subfigure}
\usepackage{verbatim}
\usepackage{wrapfig}
\usepackage{ascmac}
\usepackage{cases}
\usepackage{amsthm}
\usepackage{autobreak}

\usepackage{comment}
\usepackage {cancel}
\usepackage[normalem]{ulem}

\usepackage{microtype}

  \makeatletter
    
    \@addtoreset{equation}{section}
  \makeatother
\newtheorem{theorem}{Theorem}[section]

\newtheorem{lemma}[theorem]{Lemma}

\newcounter{c}
\newcounter{d}
\newcounter{b}

\newcommand{\ind}[1]{\mathbf{1}_{#1}}

\begin{document}
\raggedbottom
\title[Classification of fast reaction limits]{Classification of fast reaction limits for a two-component reaction--diffusion system}
\author{Yuki Tsukamoto}
\date{}

\address{Tokyo University of Science,
	162-8601, Tokyo, Japan}
\email{tsukamoto-yuki@rs.tus.ac.jp} 

\thanks{Affiliation: Tokyo University of Science,
162-8601 Tokyo, Japan.
E-mail: tsukamoto-yuki@rs.tus.ac.jp.}


\begin{abstract}
We study fast reaction limits for a two-component
reaction--diffusion system with four reaction exponents,
in which only one component diffuses.
For nonnegative segregated initial data, we classify
the limits for all exponents at least one under
homogeneous Dirichlet and Neumann boundary conditions.
The limiting equations include the heat equation,
a one-phase Stefan problem, a fixed-domain heat equation,
and a nonlinear equation with a time-integral reaction term.
We also describe the limiting behavior of the nondiffusing
component and obtain a nontrivial rescaled limit.
In the critical regime, we characterize the limit as
the greatest weak solution and establish a transition
between initial finite propagation and instantaneous positivity.
Under Dirichlet conditions, permanent zero regions
may occur for certain small initial data.
\end{abstract}
\maketitle

\section{Introduction}

Fast reaction limits arise in reaction--diffusion systems when the
reaction takes place on a much shorter time scale than diffusion.
In irreversible two-component systems, the limiting components often
become spatially segregated, and the interface between them may be
described by a Stefan-type free boundary problem.

In this paper, we consider the following system in a bounded
connected domain $\Omega\subset\mathbb R^N$ with $C^2$ boundary:
\begin{equation}\label{eq:four}
\begin{cases}
\partial_tu_k=\Delta u_k-ku_k^{m_1}v_k^{m_2},\\
\partial_tv_k=-ku_k^{m_3}v_k^{m_4},
\end{cases}
\qquad m_1,m_2,m_3,m_4\ge1.
\end{equation}
Here $k>0$ is the reaction-rate parameter.
The component $u_k$ diffuses, whereas $v_k$ remains immobile.

Various related fast-reaction systems have been studied in the literature.
Hilhorst, van der Hout and Peletier studied fast reaction limits
for systems with one diffusing and one immobile component
\cite{HilhorstEtAl1996FastReactionLimit} and subsequently considered
general monotone reaction terms
\cite{HilhorstVanDerHoutPeletier1997GeneralMonotoneFastReaction}.
Self-similar limits on unbounded domains were investigated in
\cite{CrooksHilhorst2016SelfsimilarFastreactionLimits}.
A general framework for fast reaction limits was developed in
\cite{Murakawa2021FastReactionLimit}.
More recently, Crooks and Du studied fast reaction limits with
nonlinear diffusion and the self-similar solutions and free boundaries
of the resulting limit problems
\cite{CrooksDu2024FastreactionLimitReaction,CrooksDu2025SelfSimilarSolutionsFast}.

Recent studies have also considered different reaction
structures and analytical approaches to fast reaction limits.
For nonmonotone reaction laws, oscillations can be described by
Young measures
\cite{PerthameSkrzeczkowski2023FastReactionLimit,Skrzeczkowski2022FastReactionLimit}.
For fast reversible reactions, Stephan established EDP-convergence
for a linear reaction--diffusion system using its
gradient-flow structure
\cite{Stephan2021EDPconvergenceLinearReactiondiffusion}.
Kuehn and Sulzbach developed a slow-manifold approach to
fast-reaction systems in infinite-dimensional spaces
\cite{KuehnSulzbach2025FastReactionsSlowManifolds}.
Morgan, Soresina, Tang and Tran derived a cross-diffusion limit
and convergence rates for a plant-growth model with fast reactions
\cite{MorganEtAl2026SingularLimitConvergenceRate}.

In many of these studies, the fast reaction terms in the two equations
are proportional.
Taking a suitable linear combination of the equations then
eliminates the reaction term and provides a useful starting
point for identifying the limit.
This cancellation is generally unavailable when the two
components are consumed according to different nonlinear
reaction laws.
For monomial reaction terms with four exponents,
Iida, Monobe, Murakawa and Ninomiya
\cite{IidaEtAl2017VanishingMovingImmovable}
studied representative families in which three exponents are fixed
equal to one. They obtained instantaneous interface disappearance,
finite-speed motion described by Stefan problems, and stationary
interfaces, depending on the varying exponent.
Hayashi also derived vanishing, moving, and stationary interfaces
from two-species exclusion processes with unbalanced reaction rates
in the corresponding cases with integer exponents
\cite{Hayashi2021SpatialsegregationLimitExclusion}.
Subsequent works established further results on interface
disappearance and stationary interfaces
\cite{Tsukamoto2025InterfaceDisappearanceFast,Tsukamoto2026FormationStationaryInterfaces,Tsukamoto2026VanishingInterfacesAsymmetric}.
These studies show that the relation between the nonlinearities in
the two equations can substantially change the evolution of the
limiting interface.

Fast reaction limits can also exhibit phenomena beyond
the motion of an interface.
In \cite{Tsukamoto2026PersistenceReactionTermUnbalanced}, the author
considered the exponent range
$m_2=m_4-1$ and $m_1-m_3\ge1$.
Although the nondiffusing component vanishes as $k\to\infty$,
the reaction term has a nonzero limit, and the equation for
the diffusive component retains a time-integral term.
In particular, the limiting equation is not simply the heat equation
in the region where the diffusive component is positive.
This suggests that, for a general choice of four exponents, one
should determine the limiting equation as well as the behavior of
the interface.

The purpose of the present paper is to classify the fast reaction
limit of the four-exponent system for every
$m_1,m_2,m_3,m_4\ge1$, under both homogeneous Dirichlet and Neumann
boundary conditions.
The classification depends on two differences of the
reaction exponents.
Four types of limiting equations arise: the heat equation,
a one-phase Stefan problem, a fixed-domain heat equation,
and a nonlinear equation with a time-integral reaction term.
We establish convergence of the diffusive component throughout the
exponent range and describe the corresponding behavior of the
nondiffusing component.
In a borderline case, the rescaled diffusive component
has a nontrivial limit.
Together with the limit of the nondiffusing component,
it satisfies a coupled system in the region initially
occupied by the immobile component.
For the fixed-domain heat limit, no regularity is imposed on the
boundary of that region; the limiting problem is formulated in a
variational sense.

In the critical regime $m_2=m_4-1$ and $m_1-m_3>0$,
we prove strong $L^1$ convergence of the reaction term
and characterize the limit as the greatest weak solution
of the resulting equation.
We also show that its propagation properties change at
$m_1-m_3=1$.
When $0<m_1-m_3<1$, the limiting diffusive component
remains zero almost everywhere on each compact subset
of the region initially occupied by the immobile component
for some initial time interval.
Nevertheless, it penetrates this region at arbitrarily
small positive times.
Under the Dirichlet boundary condition, we also prove
the existence of permanent zero regions for certain
sufficiently small initial data.
In contrast, instantaneous positivity holds when
$m_1-m_3\ge1$.

We now introduce the notation and limiting problems needed
to state our main results.
We set $\delta=m_1-m_3$ and $\sigma=m_2-m_4+1$.
We assume
\begin{equation}\label{eq:data}
u_0,v_0\in C(\overline\Omega),\qquad u_0,v_0\ge0,\qquad
u_0v_0=0,\qquad u_0\not\equiv0,\qquad v_0\not\equiv0.
\end{equation}
For each $k>0$, we consider \eqref{eq:four} with
initial data $(u_0,v_0)$.
We impose one of the following boundary conditions
on $u_k$, with the same choice for all $k>0$:
\begin{equation}\label{eq:bc}
\mathcal{B}u_k=0
\quad\text{on }\partial\Omega\times(0,\infty),
\qquad
\mathcal{B}u=
\begin{cases}
u & \text{Dirichlet},\\
\partial_\nu u & \text{Neumann}.
\end{cases}
\end{equation}
In the Dirichlet case, we further assume
$u_0=0$ on $\partial\Omega$.
No boundary condition is imposed on $v_k$.
By \cite[Theorem~2.1]{Tsukamoto2026PersistenceReactionTermUnbalanced},
for each $k>0$ the problem has a unique global
nonnegative solution satisfying, for every $r>N$,
\begin{equation}\label{eq:finite-regularity}
\begin{aligned}
u_k&\in C([0,\infty);C(\overline\Omega))
\cap C^1((0,\infty);C(\overline\Omega))
\cap C((0,\infty);W^{2,r}(\Omega)),\\
v_k&\in C^1([0,\infty);C(\overline\Omega)).
\end{aligned}
\end{equation}
Set
\[
D=\{x\in\Omega:v_0(x)>0\},\qquad
Q_T=\Omega\times(0,T),\qquad
R_k=ku_k^{m_1}v_k^{m_2}.
\]
For a function $w=w(x,t)$, we write $w(t)=w(\cdot,t)$.
For functions on $\Omega$, we write
$\|f\|_r=\|f\|_{L^r(\Omega)}$ for $1\le r\le\infty$.
Let $H$ be the solution of the heat equation
\[
\begin{cases}
\partial_t H=\Delta H
&\text{in }\Omega\times(0,\infty),\\
\mathcal BH=0
&\text{on }\partial\Omega\times(0,\infty),\\
H(0)=u_0
&\text{in }\Omega.
\end{cases}
\]

We put
\[
X_{\mathcal B}=
\begin{cases}
H_0^1(\Omega)&\text{Dirichlet},\\
H^1(\Omega)&\text{Neumann}.
\end{cases}
\]
We denote its dual space by $X_{\mathcal B}'$.
Since no regularity is assumed for $\partial D$,
we define the fixed-domain heat solution $U_D$
as the weak solution satisfying, for every $T>0$
and $0<\tau<T$,
\[
\begin{gathered}
U_D\in C([0,T];L^2(\Omega))
\cap L^2(0,T;X_{\mathcal B})
\cap H^1(\tau,T;L^2(\Omega)),\\
U_D(0)=u_0,\qquad
U_D=0\quad\text{a.e. on }D\times(0,T).
\end{gathered}
\]
For every $\phi\in X_{\mathcal B}$ with
$\phi=0$ a.e. in $D$, we have
\begin{equation}\label{eq:fixed}
\int_\Omega (U_D(t)-u_0)\phi\,dx
+\int_0^t\!\!\int_\Omega
\nabla U_D\cdot\nabla\phi\,dx\,ds=0
\qquad(0\le t\le T).
\end{equation}
The existence and uniqueness of $U_D$ are proved
in Lemma~\ref{lem:fixed}.
If $D\Subset\Omega$ has smooth boundary, then $U_D$
coincides with the solution of the heat equation on
$\Omega\setminus\overline D$, with zero boundary values
on $\partial D$ and the same boundary condition as $u_k$
on $\partial\Omega$, extended by zero to $D$.

For $\sigma>0$, set $z_0=v_0^\sigma/\sigma$.
When $\delta=0$, we consider the one-phase Stefan problem
\begin{equation}\label{eq:stefan}
\begin{cases}
\partial_t E=\Delta E_+
&\text{in }\Omega\times(0,\infty),\\
\mathcal BE_+=0
&\text{on }\partial\Omega\times(0,\infty),\\
E(0)=u_0-z_0
&\text{in }\Omega,
\end{cases}
\end{equation}
where $E_+=\max\{E,0\}$.
A weak solution is a function $E$ satisfying
\[
E\in L^\infty(Q_T),\qquad
E_+\in L^2(0,T;X_{\mathcal B})
\]
for every $T>0$, and
\[
\int_0^T\!\!\int_\Omega
\left(-E\partial_t\phi
+\nabla E_+\cdot\nabla\phi\right)\,dx\,dt
=\int_\Omega(u_0-z_0)\phi(0)\,dx
\]
for every $\phi\in C^1([0,T];X_{\mathcal B})$
with $\phi(T)=0$.
The uniqueness of the weak solution is proved
in Section~\ref{sec:noncritical}.

Suppose $\sigma=0$ and $\delta>0$.
For a bounded nonnegative function $w$, set
\begin{equation}\label{eq:memory}
\begin{gathered}
A_w(x,t)=\int_0^t w(x,s)^{m_3}\,ds,\\
\mathcal F[w](x,t)=
\begin{cases}
\displaystyle\frac{w(x,t)^{m_1}}{m_2 A_w(x,t)}
& x\in D,\ A_w(x,t)>0,\\[3pt]
0 & \text{otherwise}.
\end{cases}
\end{gathered}
\end{equation}
We consider the following problem:
\begin{equation}\label{eq:memory_problem}
\begin{cases}
\partial_t w=\Delta w-\mathcal F[w]
&\text{in }\Omega\times(0,\infty),\\
\mathcal Bw=0
&\text{on }\partial\Omega\times(0,\infty),\\
w(0)=u_0
&\text{in }\Omega.
\end{cases}
\end{equation}
A nonnegative function $w$ is a weak solution of
\eqref{eq:memory_problem} if, for every $T>0$
and $0<\tau<T$,
\begin{equation}\label{eq:weakclass}
\begin{gathered}
w\in L^\infty(Q_T)\cap C([0,T];L^2(\Omega))
\cap L^2(0,T;X_{\mathcal B})
\cap H^1(\tau,T;L^2(\Omega)),\\
\mathcal F[w]\in L^1(Q_T),\qquad w(0)=u_0,
\end{gathered}
\end{equation}
and
\begin{equation}\label{eq:weakform}
\int_0^T\!\!\int_\Omega
\left(-w\partial_t \varphi+\nabla w\cdot\nabla\varphi
+\mathcal F[w]\varphi\right)\,dx\,dt
=\int_\Omega u_0\varphi(0)\,dx
\end{equation}
for every $\varphi\in C^\infty(\overline\Omega\times[0,T])$
with $\varphi(T)=0$ and, in the Dirichlet case,
$\varphi=0$ on $\partial\Omega\times[0,T]$.

A weak solution $u$ is called the greatest weak solution
if $w\le u$ almost everywhere for every weak solution $w$.

We now state the main result.
\begin{theorem}\label{thm:classification}
Assume~\eqref{eq:data} and~\eqref{eq:bc}.
For every $T>0$ and $1\le r<\infty$,
the family $u_k$ converges in $C([0,T];L^r(\Omega))$
to a function $u$.
The following assertions hold for every $0<\tau<T$.

\begin{enumerate}
\item[(i)]
If $\sigma>0$ and $\delta>0$, then $u=H$ and
$v_k\to0$ in $C([\tau,T];L^r(\Omega))$.
Moreover, $R_k\to0$ in $L^1(Q_T)$.

\item[(ii)]
If $\sigma>0$ and $\delta=0$, then $u=E_+$ and
\begin{equation}\label{eq:weakv}
\frac{v_k^\sigma}{\sigma}
\overset{*}{\rightharpoonup}(-E)_+
\quad\text{in }L^\infty(Q_T).
\end{equation}

\item[(iii)]
If $\sigma>0$ and $\delta<0$,
or $\sigma=0$ and $\delta\le0$,
or $\sigma<0$, then $u=U_D$.
Furthermore,
\[
v_k\longrightarrow
\begin{cases}
v_0 & \delta<1,\\
V_* & \delta=1,\\
0 & \delta>1.
\end{cases}
\]
The first two convergences hold in
$C([0,T];L^r(\Omega))$, and the last in
$C([\tau,T];L^r(\Omega))$.
The limit $V_*$ is characterized in
Theorem~\ref{thm:second}.
Moreover, for every $K\Subset D$,
\begin{equation}\label{eq:localzero}
\sup_{K\times[0,T]}u_k\longrightarrow0.
\end{equation}

\item[(iv)]
If $\sigma=0$ and $\delta>0$, then $u$ is
the greatest weak solution of
\eqref{eq:memory_problem}.
Moreover,
\[
v_k\longrightarrow0
\quad\text{a.e. on }
(D\times(0,\infty))\cap\{A_u>0\}.
\]
If $\delta\ge1$, then $v_k\to0$
in $C([\tau,T];L^r(\Omega))$.
\end{enumerate}
\end{theorem}

We next characterize the limit $V_*$
appearing in Theorem~\ref{thm:classification} (iii).
\begin{theorem}\label{thm:second}
Suppose $\sigma<0$ and $\delta=1$, and set
$W_k=k^{1/m_3}u_k$.
Then $W_k$ is nondecreasing and $v_k$ is
nonincreasing with respect to $k$.
For every $T>0$,
\[
W_k\longrightarrow W_*,
\qquad
v_k\longrightarrow V_*
\]
locally uniformly on $D\times[0,T]$.
Moreover,
\begin{equation}\label{eq:Vstar}
V_*(x,t)=
\left(
v_0(x)^{1-m_4}
+(m_4-1)\int_0^t W_*(x,s)^{m_3}\,ds
\right)^{-1/(m_4-1)}
\end{equation}
for $x\in D$ and $t\ge0$, and
\[
0<V_*(x,t)<v_0(x)
\qquad(x\in D,\ t>0).
\]
The pair $(W_*,V_*)$ satisfies, in the weak sense
on $D\times(0,\infty)$,
\begin{equation}\label{eq:starsystem}
\begin{cases}
\partial_t W_*=\Delta W_*-W_*^{m_3+1}V_*^{m_2},\\
\partial_t V_*=-W_*^{m_3}V_*^{m_4},
\end{cases}
\end{equation}
with $W_*(0)=0$ and $V_*(0)=v_0$ in $D$.
Here $V_*$ is extended by zero to $\Omega\setminus D$.
\end{theorem}
We next study the convergence of the reaction term
and the penetration of $u$ into $D$ when
$\sigma=0$ and $\delta>0$.

\begin{theorem}\label{thm:penetration}
Suppose $\sigma=0$ and $\delta>0$, and let $u$
be the limit in Theorem~\ref{thm:classification}.
For every $T>0$,
\[
R_k\longrightarrow\mathcal F[u]
\quad\text{strongly in }L^1(Q_T).
\]
Moreover, $u$ depends on $v_0$ only through $D$.
The following assertions hold.

\begin{enumerate}
\item[(i)]
If $0<\delta<1$, set
\begin{equation}\label{eq:Cdelta}
C_\delta=
\frac{\sqrt{2m_2[N(1-\delta)+2\delta]}}
{1-\delta}
\|u_0\|_\infty^{(1-\delta)/2}.
\end{equation}
For every $B_R(x_0)\Subset D$, we have
\begin{equation}\label{eq:penetration}
u(t)=0\quad\text{a.e. in }
B_{R-C_\delta\sqrt{t}}(x_0)
\qquad
\left(0<t<\frac{R^2}{4C_\delta^2}\right).
\end{equation}
In particular, for every $K\Subset D$,
there exists $T_K>0$ such that
$u=0$ a.e. on $K\times(0,T_K)$.
Nevertheless,
\begin{equation}\label{eq:nonstationary}
\int_0^T\!\!\int_D u^{m_3}\,dx\,dt>0
\qquad(T>0).
\end{equation}

\item[(ii)]
If $\delta\ge1$, then
$u\in C(\overline\Omega\times(0,\infty))$,
$u>0$ in $\Omega\times(0,\infty)$, and
\begin{equation}\label{eq:instant}
\sup_{\overline\Omega\times[\tau,T]}
|u_k-u|\longrightarrow0
\qquad(0<\tau<T).
\end{equation}
Moreover, for every $1\le r<\infty$, the following
convergences hold in
$L^r_{\mathrm{loc}}(\Omega\times(0,\infty))$:
\begin{equation}\label{eq:rescaled}
kv_k^{m_2}\longrightarrow\frac{\ind{D}}{m_2A_u},
\qquad
k^{1/m_2}v_k\longrightarrow\ind{D}(m_2A_u)^{-1/m_2}.
\end{equation}
\end{enumerate}
\end{theorem}
When $0<\delta<1$, these results do not determine
the limit of $v_k$ on
$(D\times(0,\infty))\cap\{A_u=0\}$.

We next show that, under the Dirichlet boundary
condition, a permanent zero region may remain
when $0<\delta<1$.
\begin{theorem}\label{thm:permanent-zero}
Suppose $\sigma=0$, $0<\delta<1$, and impose
the homogeneous Dirichlet boundary condition.
Let $B_{3R}(x_0)\Subset D$ and
$0\le f\in C_c^\infty(\Omega\setminus\overline D)$
with $f\not\equiv0$.
Then there exists $\varepsilon_*>0$ such that,
for every $0<\varepsilon\le\varepsilon_*$,
the limit $u$ corresponding to $u_0=\varepsilon f$
satisfies
\[
u=0\quad\text{a.e. on }
B_{2R}(x_0)\times(0,\infty).
\]
\end{theorem}
Theorem~\ref{thm:penetration} and
Theorem~\ref{thm:permanent-zero} show that,
although the diffusive component penetrates into $D$
at arbitrarily small positive times,
a nonempty region of $D$ may remain uninvaded for all time.

The paper is organized as follows.
Section~2 establishes the estimates and compactness properties
used throughout the paper.
Section~3 treats the noncritical regimes, including the
Stefan-type limit and the rescaled limit.
Section~4 treats the critical regime.
We establish strong convergence of the reaction term
and characterize the limit as the greatest weak solution.
We also investigate initial finite propagation and,
under the Dirichlet boundary condition, permanent zero regions.

\section{Common estimates}

Unless otherwise stated, $C$ and $c$ denote positive constants
independent of $k$; their values may change from line to line.
By the comparison principle,
\begin{equation}\label{eq:bounds}
0\le u_k\le H\le\|u_0\|_\infty,\qquad
0\le v_k\le v_0\le\|v_0\|_\infty.
\end{equation}
Set
\begin{equation}\label{eq:Ak}
A_k(x,t)=\int_0^t u_k(x,s)^{m_3}\,ds.
\end{equation}
Integrating the second equation of \eqref{eq:four},
we obtain
\begin{equation}\label{eq:vformula}
v_k(x,t)=
\begin{cases}
\displaystyle
v_0(x)\left(1+(m_4-1)kv_0(x)^{m_4-1}
A_k(x,t)\right)^{-1/(m_4-1)}
& m_4>1,\\[3pt]
v_0(x)e^{-kA_k(x,t)}
& m_4=1.
\end{cases}
\end{equation}
In particular, $v_k=R_k=0$ in
$(\Omega\setminus D)\times(0,\infty)$.
By~\cite[Lemma~2.3]{Tsukamoto2026PersistenceReactionTermUnbalanced},
if $\sigma\le0$ and $0<k_1<k_2$, then
\begin{equation}\label{eq:kmonotone}
u_{k_2}\le u_{k_1},\qquad
k_1^{1/m_2}v_{k_1}\le k_2^{1/m_2}v_{k_2}.
\end{equation}

We next establish estimates independent of $k$.

\begin{lemma}\label{lem:energy}
For every $t>0$,
\begin{align}
&\int_0^t\!\!\int_\Omega R_k\,dx\,ds
\le \int_\Omega (u_0-u_k(t))\,dx
\le \|u_0\|_1,
\label{eq:mass}\\
&\frac12\|u_k(t)\|_2^2
+\int_0^t\!\!\int_\Omega
\bigl(|\nabla u_k|^2+u_kR_k\bigr)\,dx\,ds
= \frac12\|u_0\|_2^2.
\label{eq:L2}
\end{align}
Moreover,
\begin{equation}\label{eq:positiveenergy}
\begin{gathered}
\|\nabla u_k(t)\|_2^2
\le\frac{\|u_0\|_2^2}{2t}
\qquad(t>0),\\
\int_\tau^T\|(u_k)_t\|_2^2\,dt
\le\frac{\|u_0\|_2^2}{4\tau}
\qquad(0<\tau<T).
\end{gathered}
\end{equation}
\end{lemma}
\begin{proof}
The following calculations are justified by
\eqref{eq:finite-regularity}.
Integrating the equation for $u_k$ in \eqref{eq:four}
over $\Omega$ and using \eqref{eq:bc} and $u_k\ge0$,
we obtain
\[
\frac{d}{dt}\int_\Omega u_k\,dx
\le-\int_\Omega R_k\,dx.
\]
Integrating over $(s,t)$ and letting $s\downarrow0$
gives \eqref{eq:mass}.
Similarly, multiplying the first equation in \eqref{eq:four}
by $u_k$, integrating over $\Omega\times(s,t)$,
and letting $s\downarrow0$ yields \eqref{eq:L2}.

To estimate $\partial_tu_k$, we multiply the first equation
in \eqref{eq:four} by $\partial_tu_k$ and integrate over
$\Omega$. Then
\[
\|\partial_tu_k(t)\|_2^2
+\frac12\frac{d}{dt}\|\nabla u_k(t)\|_2^2
=-\int_\Omega R_k\partial_tu_k\,dx.
\]
Set
\[
\mathcal E_k(t)
=\frac12\|\nabla u_k(t)\|_2^2
+\frac{k}{m_1+1}
\int_\Omega u_k(t)^{m_1+1}v_k(t)^{m_2}\,dx.
\]
Since
\[
\begin{aligned}
\frac{d}{dt}\left(
\frac{k}{m_1+1}
\int_\Omega u_k^{m_1+1}v_k^{m_2}\,dx
\right)
={}&\int_\Omega R_k\partial_tu_k\,dx\\
&+\frac{km_2}{m_1+1}
\int_\Omega u_k^{m_1+1}v_k^{m_2-1}
\partial_tv_k\,dx,
\end{aligned}
\]
the second equation in \eqref{eq:four} gives
\begin{equation}\label{eq:dissipation}
\mathcal E_k'(t)
=-\|\partial_tu_k(t)\|_2^2
-\frac{m_2k^2}{m_1+1}
\int_\Omega
u_k^{m_1+m_3+1}v_k^{m_2+m_4-1}\,dx
\le0.
\end{equation}

Since $m_1\ge1$, \eqref{eq:L2} implies
\[
\int_0^t\mathcal E_k(s)\,ds
\le\frac12\int_0^t\!\!\int_\Omega
\bigl(|\nabla u_k|^2+u_kR_k\bigr)\,dx\,ds
\le\frac14\|u_0\|_2^2.
\]
As $\mathcal E_k$ is nonincreasing,
\[
t\mathcal E_k(t)
\le\int_0^t\mathcal E_k(s)\,ds
\le\frac14\|u_0\|_2^2.
\]
Consequently,
\[
\|\nabla u_k(t)\|_2^2
\le2\mathcal E_k(t)
\le\frac{\|u_0\|_2^2}{2t}.
\]
Finally, integrating \eqref{eq:dissipation}
from $\tau$ to $T$, we obtain
\[
\int_\tau^T\|\partial_tu_k(t)\|_2^2\,dt
\le\mathcal E_k(\tau)
\le\frac{\|u_0\|_2^2}{4\tau}.
\]
This completes the proof.
\end{proof}
We also need uniform estimates near the initial time.

\begin{lemma}\label{lem:compact}
We have
\begin{equation}\label{eq:trace}
\lim_{t\downarrow0}\sup_{k>0}\|u_k(t)-u_0\|_2=0,
\qquad
\lim_{t\downarrow0}\sup_{k>0}
\int_0^t\!\!\int_\Omega R_k\,dx\,ds=0.
\end{equation}
Moreover, for every $0<\tau<T$ and
$\tau\le s<t\le T$,
\begin{equation}\label{eq:time-continuity}
\|u_k(t)-u_k(s)\|_2
\le\frac{\|u_0\|_2}{2\sqrt{\tau}}|t-s|^{1/2}.
\end{equation}
\end{lemma}
\begin{proof}
Since $u_0\ge0$ is continuous, it can be approximated
in $L^1(\Omega)$ by functions
$\phi\in C_c^\infty(\{u_0>0\})$.
Since $v_0=0$ on the support of $\phi$,
we have $R_k=0$ there.
Multiplying the first equation in \eqref{eq:four}
by $\phi$ and integrating, we obtain
\[
\left|\int_\Omega (u_k(t)-u_0)\phi\,dx\right|
\le t\|u_0\|_\infty\|\Delta\phi\|_1.
\]
By \eqref{eq:bounds} and \eqref{eq:L2},
\[
\begin{aligned}
\|u_k(t)-u_0\|_2^2
&=\|u_k(t)\|_2^2+\|u_0\|_2^2
-2\int_\Omega u_k(t)u_0\,dx\\
&\le 2\int_\Omega (u_0-u_k(t))u_0\,dx\\
&=2\int_\Omega (u_0-u_k(t))(u_0-\phi)\,dx+2\int_\Omega (u_0-u_k(t))\phi\,dx\\
&\le 2\|u_0\|_\infty\|u_0-\phi\|_1
+2t\|u_0\|_\infty\|\Delta\phi\|_1.
\end{aligned}
\]
Letting $t\downarrow0$ and then
$\phi\to u_0$ in $L^1(\Omega)$
gives the first limit in \eqref{eq:trace}.
Moreover, by \eqref{eq:mass},
\[
\begin{aligned}
\int_0^t\!\!\int_\Omega R_k\,dx\,ds
&\le \int_\Omega (u_0-u_k(t))\,dx\\
&\le |\Omega|^{1/2}\|u_k(t)-u_0\|_2.
\end{aligned}
\]
The second limit in \eqref{eq:trace} follows
from the first.

Finally, for $0<\tau\le s<t\le T$,
\eqref{eq:positiveenergy} and the
Cauchy--Schwarz inequality give
\[
\begin{aligned}
\|u_k(t)-u_k(s)\|_2
&\le \int_s^t\|\partial_tu_k(r)\|_2\,dr\\
&\le |t-s|^{1/2}
\left(\int_\tau^T
\|\partial_tu_k(r)\|_2^2\,dr\right)^{1/2}\\
&\le \frac{\|u_0\|_2}{2\sqrt{\tau}}
|t-s|^{1/2}.
\end{aligned}
\]
This completes the proof.
\end{proof}
We next identify limits that vanish on $D$.

\begin{lemma}\label{lem:fixed}
The fixed-domain heat problem \eqref{eq:fixed}
has a unique solution $U_D$.
Suppose that a subsequence of $u_k$ converges
to $U$ in $C([0,T];L^2(\Omega))$ as $k\to\infty$,
and that $U=0$ a.e. on $D\times(0,T)$.
Then $U=U_D$.
Moreover, the same subsequence converges to zero
uniformly on $K\times[0,T]$ for every $K\Subset D$.
\end{lemma}
\begin{proof}
Let
\[
X_D=\{\phi\in X_{\mathcal B}:
\phi=0\text{ a.e. in }D\}.
\]
Then $X_D$ is a closed subspace of $X_{\mathcal B}$.
The bilinear form
\[
a(\phi,\psi)=
\int_\Omega\nabla\phi\cdot\nabla\psi\,dx
\]
is bounded on $X_D$ and satisfies the coercivity estimate
\[
a(\phi,\phi)+\|\phi\|_2^2
=\|\phi\|_{H^1(\Omega)}^2.
\]
Since $u_0$ can be approximated in $L^2(\Omega)$
by functions in $C_c^\infty(\{u_0>0\})\subset X_D$,
we can apply the standard variational existence theorem
\cite{LionsMagenes1972NonHomogeneousBoundaryValue}
to $X_D$ and its closure in $L^2(\Omega)$.
This gives a solution $U_D$ of \eqref{eq:fixed}.
Since the form is closed, symmetric, and nonnegative,
its associated operator is nonnegative and self-adjoint.
The corresponding semigroup is therefore analytic
\cite[Corollary~II.4.8]{EngelNagel2006ShortCourseOperatorSemigroups}.
Consequently,
\[
U_D\in H^1(\tau,T;L^2(\Omega))
\qquad(0<\tau<T).
\]
The maximum principle also gives
$0\le U_D\le\|u_0\|_\infty$.

To prove uniqueness, let $U_1$ and $U_2$ be two
solutions and set $W=U_1-U_2$.
Subtracting their equations, multiplying
the resulting equation by $W$, and integrating
over $(s,t)$, we obtain
\[
\frac12\|W(t)\|_2^2
+\int_s^t\|\nabla W(r)\|_2^2\,dr
=\frac12\|W(s)\|_2^2.
\]
Letting $s\downarrow0$ and using $W(0)=0$,
we obtain $W=0$.

Next, suppose that a subsequence of $u_k$
converges to $U$ in $C([0,T];L^2(\Omega))$
and that $U=0$ a.e. on $D\times(0,T)$.
We retain the notation $u_k$ for this subsequence.
By \eqref{eq:L2}, \eqref{eq:positiveenergy}
and weak compactness, we obtain
\[
U\in L^2(0,T;X_D)
\cap H^1(\tau,T;L^2(\Omega))
\qquad(0<\tau<T).
\]
Moreover, $U(0)=u_0$.
Since $R_k=0$ outside $D$, we have
$R_k\phi=0$ for every $\phi\in X_D$.
Multiplying the first equation in \eqref{eq:four}
by $\phi$ and integrating, we obtain
\[
\int_\Omega(u_k(t)-u_0)\phi\,dx
+\int_0^t\!\!\int_\Omega
\nabla u_k\cdot\nabla\phi\,dx\,ds=0.
\]
Passing to the limit gives \eqref{eq:fixed} for $U$,
and hence $U=U_D$ by uniqueness.

Finally, since $R_k\ge0$, $u_k$ is a nonnegative
subsolution of the heat equation.
Let $K\Subset G\Subset D$ and $0<\tau<T$.
The interior parabolic estimate
\cite{LadyzhenskayaSolonnikovUraltseva1968LinearQuasilinearParabolic}
gives
\[
\sup_{K\times[\tau,T]}u_k
\le C\int_{\tau/2}^T\!\!\int_G u_k\,dx\,dt
\longrightarrow0,
\]
where we used $u_k\to U$ in
$C([0,T];L^2(\Omega))$ and
$U=0$ a.e. on $D\times(0,T)$.
On the other hand, $u_k\le H$ and
$H(t)\to u_0=0$ uniformly on $K$ as $t\downarrow0$.
Therefore,
\[
\limsup_{k\to\infty}
\sup_{K\times[0,T]}u_k
\le \sup_{K\times[0,\tau]}H
\longrightarrow0
\qquad\text{as }\tau\downarrow0.
\]
This completes the proof.
\end{proof}

\section{Noncritical regimes}\label{sec:noncritical}

We first consider the case $\sigma>0$. Set
\[
z_k=\frac{v_k^\sigma}{\sigma},\qquad
\beta=\frac{m_2}{\sigma}\ge1.
\]
The equations in \eqref{eq:four} give
\[
\begin{aligned}
\partial_t u_k&=\Delta u_k-R_k,\\
\partial_t z_k&=-k\sigma^\beta u_k^{m_3}z_k^\beta,\\
R_k&=k\sigma^\beta u_k^{m_1}z_k^\beta.
\end{aligned}
\]
In particular,
\begin{equation}\label{eq:zmass}
\int_0^T\!\!\int_\Omega(-\partial_tz_k)\,dx\,dt
=\int_\Omega(z_0-z_k(T))\,dx
\le\|z_0\|_1.
\end{equation}

We begin with the following estimate.

\begin{lemma}\label{lem:cutoff}
For every $T>0$ and $\eta>0$,
\begin{equation}\label{eq:cutoff}
\int_{Q_T\cap\{u_k\ge\eta\}}
(-\partial_tz_k)\,dx\,dt\longrightarrow0
\qquad\text{as }k\to\infty.
\end{equation}
\end{lemma}
\begin{proof}
Let $\chi_\eta$ be a nondecreasing smooth function
such that $\chi_\eta=0$ on $[0,\eta/2]$,
$\chi_\eta=1$ on $[\eta,\infty)$, and
$0\le\chi_\eta'\le4/\eta$.
Since $u_0=0$ on $D$ and $H(t)\to u_0$
uniformly as $t\downarrow0$, \eqref{eq:bounds}
implies that there exists $t_\eta\in(0,T)$,
independent of $k$, such that
\[
u_k(x,t)\le H(x,t)<\eta/2
\qquad(x\in D,\ 0\le t\le t_\eta).
\]
Since $z_k=0$ outside $D$, it follows that
the integrand in \eqref{eq:cutoff} vanishes
for $t\le t_\eta$ and
$\chi_\eta(u_k(t_\eta))z_k(t_\eta)=0$.
Integrating by parts in time, we obtain
\begin{equation}\label{eq:cutoffibp}
\begin{aligned}
&\int_{Q_T\cap\{u_k\ge\eta\}}
(-\partial_tz_k)\,dx\,dt\\
\le&
\int_{t_\eta}^T\!\!\int_\Omega
\chi_\eta(u_k)(-\partial_tz_k)\,dx\,dt\\
=&
-\int_\Omega\chi_\eta(u_k(T))z_k(T)\,dx
+\int_{t_\eta}^T\!\!\int_\Omega
\chi_\eta'(u_k)(u_k)_tz_k\,dx\,dt\\
\le&\frac4\eta
\|(u_k)_t\|_{L^2(\Omega\times(t_\eta,T))}
\|z_k\ind{\{u_k>\eta/2\}}\|_{L^2(Q_T)}.
\end{aligned}
\end{equation}
By \eqref{eq:positiveenergy}, the first factor
on the right-hand side is bounded independently of $k$.

On the other hand, since
$-\partial_tz_k=k\sigma^\beta u_k^{m_3}z_k^\beta$,
\eqref{eq:zmass} gives
\[
\int_{Q_T\cap\{u_k>\eta/2\}}z_k^\beta\,dx\,dt
\le\frac{2^{m_3}}{k\sigma^\beta\eta^{m_3}}
\int_0^T\!\!\int_\Omega(-\partial_tz_k)\,dx\,dt
\le\frac{2^{m_3}\|z_0\|_1}
{k\sigma^\beta\eta^{m_3}}.
\]
If $1\le\beta\le2$, then $0\le z_k\le z_0$ implies
\[
\int_{Q_T\cap\{u_k>\eta/2\}}z_k^2\,dx\,dt
\le\|z_0\|_\infty^{2-\beta}
\int_{Q_T\cap\{u_k>\eta/2\}}z_k^\beta\,dx\,dt
\longrightarrow0.
\]
If $\beta>2$, H\"older's inequality gives
\[
\int_{Q_T\cap\{u_k>\eta/2\}}z_k^2\,dx\,dt
\le |Q_T|^{1-2/\beta}
\left(\int_{Q_T\cap\{u_k>\eta/2\}}
z_k^\beta\,dx\,dt\right)^{2/\beta}
\longrightarrow0.
\]
Thus, in both cases,
\[
\|z_k\ind{\{u_k>\eta/2\}}\|_{L^2(Q_T)}
\longrightarrow0.
\]
Combining this with \eqref{eq:cutoffibp}
gives \eqref{eq:cutoff}.
This completes the proof.
\end{proof}

\begin{proof}[Proof of Theorem~\ref{thm:classification} for $\sigma>0$]
For each $t>0$, \eqref{eq:bounds} and
\eqref{eq:positiveenergy} imply that $\{u_k(t)\}$
is bounded in $H^1(\Omega)$, and hence relatively
compact in $L^2(\Omega)$ by the Rellich--Kondrachov
theorem.
Moreover, Lemma~\ref{lem:compact} gives
uniform equicontinuity in time, including at $t=0$.
Since $u_k(0)=u_0$, the Arzel\`a--Ascoli theorem
implies that $\{u_k\}$ is relatively compact in
$C([0,T];L^2(\Omega))$.
By \eqref{eq:bounds}, this relative compactness
also holds in $C([0,T];L^r(\Omega))$
for every finite $r$.

Suppose $\delta>0$. Since
$R_k=-u_k^\delta\partial_tz_k$,
\eqref{eq:zmass} gives
\[
\int_{Q_T}R_k\,dx\,dt
\le\eta^\delta\|z_0\|_1
+\|u_0\|_\infty^\delta
\int_{Q_T\cap\{u_k\ge\eta\}}
(-\partial_tz_k)\,dx\,dt.
\]
By Lemma~\ref{lem:cutoff}, letting $k\to\infty$
and then $\eta\downarrow0$ yields
$R_k\to0$ in $L^1(Q_T)$.
Since $H-u_k\ge0$ has zero initial and boundary data
and satisfies
\[
\partial_t(H-u_k)=\Delta(H-u_k)+R_k,
\]
the $L^1$ estimate for the heat equation gives
\begin{equation}\label{eq:heatconvergence}
\sup_{0\le t\le T}\|H(t)-u_k(t)\|_1
\le\|R_k\|_{L^1(Q_T)}
\longrightarrow0.
\end{equation}
By \eqref{eq:bounds}, this convergence also holds in
$C([0,T];L^r(\Omega))$ for every finite $r$.

Fix $0<\tau<T$. By \eqref{eq:heatconvergence},
\[
A_k(\tau)\longrightarrow
\int_0^\tau H(s)^{m_3}\,ds
\qquad\text{in }L^1(\Omega).
\]
The strong maximum principle implies that the limiting
function is positive throughout $\Omega$.
By \eqref{eq:vformula}, we obtain
$v_k(\tau)\to0$ in $L^1(\Omega)$.
By \eqref{eq:bounds}, this convergence also holds
in $L^r(\Omega)$ for every finite $r$.
Since $v_k$ is nonincreasing in time,
\[
\sup_{\tau\le t\le T}\|v_k(t)\|_r
\le\|v_k(\tau)\|_r\longrightarrow0.
\]

Suppose $\delta=0$. Then $m_1=m_3$ and
$R_k=-\partial_tz_k$, so
\begin{equation}\label{eq:enthalpyapprox}
\partial_t(u_k-z_k)=\Delta u_k.
\end{equation}
After passing to a subsequence, we have
$u_k\to U$ in $C([0,T];L^r(\Omega))$
for every finite $r$, and
$z_k\overset{*}{\rightharpoonup}Z$
in $L^\infty(Q_T)$.
Both limits are nonnegative.
By \eqref{eq:zmass},
\[
\int_{Q_T}u_k^{m_3}z_k^\beta\,dx\,dt
\le\frac{\|z_0\|_1}{k\sigma^\beta}.
\]
Since $u_k$ and $z_k$ are uniformly bounded,
\[
(u_kz_k)^{\max\{m_3,\beta\}}
\le C u_k^{m_3}z_k^\beta.
\]
Consequently, $u_kz_k\to0$ in $L^1(Q_T)$.
The strong convergence of $u_k$ and weak-* convergence
of $z_k$ give $UZ=0$ a.e. in $Q_T$.
Setting $E=U-Z$, we obtain $U=E_+$ and $Z=(-E)_+$.
By \eqref{eq:L2}, we also have
$E_+=U\in L^2(0,T;X_{\mathcal B})$.
Together with \eqref{eq:enthalpyapprox}, the above
convergence shows that $E$ is a weak solution of
\eqref{eq:stefan} with initial datum $u_0-z_0$.

To prove uniqueness, let $E_1,E_2$ be two bounded
weak solutions, and set
\[
\xi=E_1-E_2,\qquad
\zeta=(E_1)_+-(E_2)_+.
\]
Then $\partial_t \xi=\Delta\zeta$ and $\xi(0)=0$.
In the Neumann case, integration over $\Omega$
gives
\[
\int_\Omega \xi(x,t)\,dx=0.
\]
Let $G_{\mathcal B}\xi$ solve
\[
-\Delta(G_{\mathcal B}\xi)=\xi
\quad\text{in }\Omega,
\]
with
\[
\begin{cases}
G_{\mathcal B}\xi=0
&\text{on }\partial\Omega
\quad\text{(Dirichlet)},\\
\partial_\nu(G_{\mathcal B}\xi)=0
&\text{on }\partial\Omega
\quad\text{(Neumann)}.
\end{cases}
\]
In the Neumann case, we also impose
$\int_\Omega G_{\mathcal B}\xi\,dx=0$.

Since $\zeta\in L^2(0,T;X_{\mathcal B})$,
we have $\partial_t \xi \in L^2(0,T;X_{\mathcal B}')$.
By the symmetry of $G_{\mathcal B}$ and
integration by parts, we obtain
\[
\begin{aligned}
\frac12\frac{d}{dt}
\int_\Omega \xi\,G_{\mathcal B}\xi\,dx
&=\langle\partial_t \xi,G_{\mathcal B}\xi\rangle\\
&=-\int_\Omega
\nabla\zeta\cdot\nabla(G_{\mathcal B}\xi)\,dx\\
&=-\int_\Omega\xi\zeta\,dx\le0,
\end{aligned}
\]
where the last inequality follows from
$\xi\zeta\ge0$.
Moreover,
\[
\int_\Omega \xi\,G_{\mathcal B}\xi\,dx
=\int_\Omega|\nabla(G_{\mathcal B}\xi)|^2\,dx
\ge0.
\]
Integrating in time and using $\xi(0)=0$, we obtain
\[
0\le\int_\Omega|\nabla(G_{\mathcal B}\xi(t))|^2\,dx
\le\int_\Omega|\nabla(G_{\mathcal B}\xi(0))|^2\,dx
=0.
\]
Hence $G_{\mathcal B}\xi=0$ by the Dirichlet
boundary condition or the zero-mean condition
in the Neumann case.
Thus $\xi=0$ and $E_1=E_2$.
Since the limit is unique, $u_k\to E_+$
for the whole family.
The weak-* convergence \eqref{eq:weakv}
follows from $Z=(-E)_+$.

Finally, suppose $\delta<0$. Since
$-\partial_tz_k=u_k^{-\delta}R_k$,
\eqref{eq:mass} gives
\[
\int_{Q_T}(-\partial_tz_k)\,dx\,dt
\le\eta^{-\delta}\|u_0\|_1
+\int_{Q_T\cap\{u_k\ge\eta\}}
(-\partial_tz_k)\,dx\,dt.
\]
By Lemma~\ref{lem:cutoff}, letting $k\to\infty$
and then $\eta\downarrow0$ yields
$-\partial_tz_k\to0$ in $L^1(Q_T)$.
Since $z_k$ is nonincreasing in time,
\[
\sup_{0\le t\le T}\|z_0-z_k(t)\|_1
=\int_{Q_T}(-\partial_tz_k)\,dx\,dt
\longrightarrow0.
\]
By boundedness, $z_k\to z_0$ in
$C([0,T];L^r(\Omega))$ for every finite $r$.
Hence $v_k\to v_0$ in the same space.

Let $u_k\to U$ in $C([0,T];L^2(\Omega))$
for a subsequence. By \eqref{eq:zmass},
\[
\int_{Q_T}u_k^{m_3}z_k^\beta\,dx\,dt
\le\frac{\|z_0\|_1}{k\sigma^\beta}.
\]
The strong convergence of $u_k$ and $z_k$,
together with their uniform boundedness, gives
$U^{m_3}z_0^\beta=0$ a.e. in $Q_T$.
Since $z_0>0$ on $D$, we have $U=0$
a.e. on $D\times(0,T)$.
Lemma~\ref{lem:fixed} gives $U=U_D$
and \eqref{eq:localzero} for this subsequence.
Since the limit is unique, the convergence
holds for the whole family. This completes the proof.
\end{proof}

We now consider the case $\sigma<0$. Set
\[
q=m_4-1>m_2,\qquad
\lambda=1-\frac{m_2}{q}>0.
\]
Fix $G\Subset D$ and $T>0$.
By \eqref{eq:vformula} and
$A_k\le T\|u_0\|_\infty^{m_3}$,
there exists $c_{G,T}>0$, independent of $k$,
such that
\begin{equation}\label{eq:negativeabsorption}
v_k(x,t)\ge c_{G,T}k^{-1/q},\qquad
R_k(x,t)\ge c_{G,T}k^\lambda u_k(x,t)^{m_1}
\end{equation}
for all $(x,t)\in G\times[0,T]$ and $k\ge1$.
By \eqref{eq:mass},
$u_k\to0$ in $L^{m_1}(G\times(0,T))$.
Thus every limit $U$ of $u_k$ satisfies
$U=0$ a.e. on $D\times(0,T)$.
The relative compactness established above and
Lemma~\ref{lem:fixed} therefore give
$u_k\to U_D$ in $C([0,T];L^r(\Omega))$
for every finite $r$, together with
\eqref{eq:localzero}.
To determine the limit of $v_k$, we use the following estimate.

\begin{lemma}\label{lem:absorption}
Let $B_R(x_0)\subset\mathbb R^N$, $T>0$,
and $\Lambda>0$.
Suppose that $w\ge0$ is a bounded weak subsolution of
\[
\partial_t w \le\Delta w-\Lambda w^{m_1}
\quad\text{in }B_R(x_0)\times(0,T),
\]
satisfying
\[
w(t)\to0\quad\text{in }L^2(B_R(x_0))
\quad\text{as }t\downarrow0,
\]
and
\[
w\in L^2(s,T;H^1(B_R(x_0)))
\cap H^1(s,T;L^2(B_R(x_0)))
\qquad(0<s<T).
\]
If $m_1>1$, then
\begin{equation}\label{eq:absorptionpower}
\operatorname*{ess\,sup}_{B_{R/2}(x_0)\times(0,T)}w
\le C(\Lambda R^2)^{-1/(m_1-1)}.
\end{equation}
If $m_1=1$, then
\begin{equation}\label{eq:absorptionlinear}
\operatorname*{ess\,sup}_{B_{R/2}(x_0)\times(0,T)}w
\le C\|w\|_\infty e^{-cR\sqrt\Lambda}.
\end{equation}
Here $\|w\|_\infty$ denotes the
$L^\infty(B_R(x_0)\times(0,T))$ norm,
and $C,c>0$ depend only on $m_1$ and $N$.
\end{lemma}
\begin{proof}
By translation, take $x_0=0$.
Suppose $m_1>1$.
Set
\[
\rho=\frac{2}{m_1-1},\qquad
C_0=\bigl(2N\rho+4\rho(\rho+1)\bigr)^{1/(m_1-1)}
\]
and
\[
\Phi(x)=C_0(\Lambda R^2)^{-1/(m_1-1)}
\left(1-\frac{|x|^2}{R^2}\right)^{-\rho}.
\]
Since $\rho m_1=\rho+2$, direct calculation gives
\[
\begin{aligned}
\Delta\Phi
&\le
\frac{C_0(\Lambda R^2)^{-1/(m_1-1)}}{R^2}
[2N\rho+4\rho(\rho+1)]
\left(1-\frac{|x|^2}{R^2}\right)^{-\rho-2},\\
\Lambda\Phi^{m_1}
&=
\frac{C_0^{m_1}(\Lambda R^2)^{-1/(m_1-1)}}{R^2}
\left(1-\frac{|x|^2}{R^2}\right)^{-\rho-2}.
\end{aligned}
\]
By the choice of $C_0$, we have
$\Delta\Phi\le\Lambda\Phi^{m_1}$.
Since $\Phi$ diverges as $|x|\uparrow R$,
it dominates $w$ on the boundaries of
concentric balls with radii sufficiently close to $R$.
Applying the comparison principle on these balls
and letting their radii tend to $R$, we obtain
$w\le\Phi$ a.e. in $B_R\times(0,T)$.
Since $1-|x|^2/R^2\ge3/4$ on $B_{R/2}$,
we obtain \eqref{eq:absorptionpower}.

Suppose $m_1=1$.
Set $\alpha=\sqrt{\Lambda/N}$ and
\[
\Psi(x)=\|w\|_\infty
\frac{\cosh(\alpha|x|)}{\cosh(\alpha R)}.
\]
Since $\tanh s\le s$ for $s\ge0$,
we have
\[
\begin{aligned}
\Delta\Psi
&=\left(
\alpha^2
+(N-1)\frac{\alpha}{|x|}
\tanh(\alpha|x|)
\right)\Psi\\
&\le N\alpha^2\Psi
=\Lambda\Psi
\end{aligned}
\]
for $x\ne0$.
The same inequality holds at $x=0$ by continuity.
Moreover, $\Psi=\|w\|_\infty$ on $\partial B_R$.
The comparison principle and the zero initial trace
give $w\le\Psi$ a.e. in $B_R\times(0,T)$.
Since
\[
\frac{\cosh(\alpha R/2)}{\cosh(\alpha R)}
\le 2e^{-\alpha R/2},
\]
we obtain
\[
\operatorname*{ess\,sup}_{B_{R/2}\times(0,T)}w
\le 2\|w\|_\infty
e^{-R\sqrt{\Lambda}/(2\sqrt N)}.
\]
This proves \eqref{eq:absorptionlinear}
and completes the proof.
\end{proof}

\begin{proof}[Proof of Theorem~\ref{thm:classification}(iii)
for $\sigma<0$ and $\delta\ne1$]
Suppose $\delta<1$.
Fix $K\Subset D$ and $T>0$.
If $m_1>1$, then
$0<m_1-1=m_3+\delta-1<m_3$.
By \eqref{eq:negativeabsorption} and
Lemma~\ref{lem:absorption}, we have
\[
u_k\le Ck^{-\alpha_0},\qquad
\alpha_0=\frac{\lambda}{m_1-1},
\]
on a compact neighborhood of $K$ in $D$.

Suppose that $u_k\le C_jk^{-\alpha_j}$
on a compact neighborhood of $K$
for $0\le t\le T$.
By \eqref{eq:vformula}, we obtain
\[
v_k\ge c_jk^{-\beta_j},\qquad
\beta_j=\frac{(1-m_3\alpha_j)_+}{q}
\]
on the same neighborhood.
Thus $R_k\ge c_j^{m_2}k^{1-m_2\beta_j}u_k^{m_1}$.
Applying Lemma~\ref{lem:absorption} again gives
\[
u_k\le C_{j+1}k^{-\alpha_{j+1}},\qquad
\alpha_{j+1}=\frac{1-m_2\beta_j}{m_1-1}
\]
on a smaller neighborhood of $K$.

As long as $\alpha_j\le1/m_3$, we have
\[
\begin{aligned}
\alpha_{j+1}-\alpha_j
&=\frac{\lambda+
\bigl(m_3(1-\lambda)-(m_1-1)\bigr)\alpha_j}{m_1-1}\\
&\ge
\frac{\min\{\lambda,\,1-(m_1-1)/m_3\}}{m_1-1}
>0.
\end{aligned}
\]
Thus, after finitely many steps, we obtain
$u_k\le Ck^{-\alpha_j}$ with $\alpha_j>1/m_3$
on a neighborhood of $K$.
By \eqref{eq:vformula}, $v_k$ is then bounded below
by a positive constant, independent of $k$,
on this neighborhood.
Applying Lemma~\ref{lem:absorption} once more, we obtain
\[
u_k\le Ck^{-1/(m_1-1)}
\quad\text{on }K\times[0,T].
\]
Integrating the second equation of \eqref{eq:four},
we obtain
\[
\begin{aligned}
0\le v_0(x)-v_k(x,t)
&\le k\|v_0\|_\infty^{m_4}
\int_0^T u_k(x,s)^{m_3}\,ds\\
&\le Ck^{1-m_3/(m_1-1)}
\longrightarrow0,
\end{aligned}
\]
uniformly on $K\times[0,T]$, since $m_1-1<m_3$.

If $m_1=1$, choose $G$ with
$K\Subset G\Subset D$.
By \eqref{eq:negativeabsorption} and
\eqref{eq:absorptionlinear} with
$\Lambda=c_{G,T}k^\lambda$, we obtain
\begin{equation}\label{eq:lineardecay}
u_k\le Ce^{-ck^{\lambda/2}}
\quad\text{on }K\times[0,T].
\end{equation}
Integrating the second equation of \eqref{eq:four},
we obtain
\[
\begin{aligned}
0\le v_0(x)-v_k(x,t)
&\le k\|v_0\|_\infty^{m_4}
\int_0^T u_k(x,s)^{m_3}\,ds\\
&\le Ck e^{-ck^{\lambda/2}}
\longrightarrow0
\end{aligned}
\]
uniformly on $K\times[0,T]$.
Since $0\le v_0-v_k\le v_0$ and
$v_k=v_0=0$ outside $D$,
exhaustion by compact subsets of $D$ gives
$v_k(T)\to v_0$ in $L^r(\Omega)$
for every finite $r$.
Since $v_k$ is nonincreasing in time,
we obtain $v_k\to v_0$ in
$C([0,T];L^r(\Omega))$.

Suppose $\delta>1$.
Then $m_1-1=m_3+\delta-1>m_3$.
Let $w$ solve
\begin{equation}\label{eq:scalar}
\begin{cases}
\partial_t w=\Delta w-\|v_0\|_\infty^{m_2}w^{m_1},\\
\mathcal Bw=0,\qquad w(0)=u_0.
\end{cases}
\end{equation}
By the comparison and strong maximum principles,
we have
$0\le w\le\|u_0\|_\infty$ and
$w(x,t)>0$ for $x\in\Omega$ and $t>0$.

For $k\ge1$, set
$\underline u_k=k^{-1/(m_1-1)}w$.
Then
\[
\partial_t \underline u_k-\Delta\underline u_k
+k\|v_0\|_\infty^{m_2}\underline u_k^{m_1}=0,
\qquad
\underline u_k(0)\le u_0.
\]
Since $v_k\le\|v_0\|_\infty$,
the comparison principle gives
\begin{equation}\label{eq:scalarlower}
u_k\ge k^{-1/(m_1-1)}w.
\end{equation}
For every $K\Subset\Omega$ and $0<\tau<T$,
\[
\begin{aligned}
kA_k(x,t)
&\ge k^{1-m_3/(m_1-1)}
\int_0^t w(x,s)^{m_3}\,ds\\
&\ge c k^{1-m_3/(m_1-1)}
\longrightarrow\infty
\end{aligned}
\]
uniformly on $K\times[\tau,T]$.
By \eqref{eq:vformula},
$v_k\le(qkA_k)^{-1/q}$ on $D$,
while $v_k=0$ outside $D$.
Hence $v_k\to0$ locally uniformly on
$\Omega\times(0,\infty)$.
Since $0\le v_k\le v_0$,
exhaustion by compact subsets of $D$ gives
$v_k\to0$ in $C([\tau,T];L^r(\Omega))$
for every finite $r$. This
completes the proof.
\end{proof}

\begin{proof}[Proof of Theorem~\ref{thm:second}]
Since $m_1=m_3+1$, the rescaling
$W_k=k^{1/m_3}u_k$ gives
\begin{equation}\label{eq:scaled}
\begin{cases}
\partial_t W_k=\Delta W_k-W_k^{m_3+1}v_k^{m_2},\\
\partial_t v_k=-W_k^{m_3}v_k^{q+1},\\
W_k(0)=k^{1/m_3}u_0,\quad v_k(0)=v_0,\quad
\mathcal BW_k=0.
\end{cases}
\end{equation}
Since only the first initial datum depends on $k$,
the comparison principle
\cite[Theorem~2.2]{Tsukamoto2026PersistenceReactionTermUnbalanced}
gives
\[
W_{k_1}\le W_{k_2},\qquad
v_{k_2}\le v_{k_1}
\qquad(0<k_1<k_2).
\]

Fix $T>0$ and $B_R(x_0)\Subset D$.
Choose $0<b_0\le\min_{\overline{B_R(x_0)}}v_0$.
By translation, take $x_0=0$.
Set
\[
\rho=\frac{2}{m_3\lambda},\qquad
\Phi(x)=C_0(R^2-|x|^2)^{-\rho},
\]
and
\[
Z(x,t)=
\left[b_0^{-q}+qt\Phi(x)^{m_3}\right]^{-1/q}.
\]
Choose $C_0>0$ large enough that
$\Phi\ge1$ in $B_R$.
For $0\le t\le T$,
\[
\begin{aligned}
\Phi^{m_3+1}Z^{m_2}
&\ge(b_0^{-q}+qT)^{-(1-\lambda)}
\Phi^{1+m_3\lambda},\\
\Delta\Phi
&\le C_0R^2[2N\rho+4\rho(\rho+1)]
(R^2-|x|^2)^{-\rho-2}.
\end{aligned}
\]
Since $\rho(1+m_3\lambda)=\rho+2$,
increasing $C_0$ if necessary gives
$\Delta\Phi\le\Phi^{m_3+1}Z^{m_2}$.
Moreover, the definition of $Z$ gives
$\partial_t Z=-\Phi^{m_3}Z^{q+1}$.
Thus $\Phi$ is an upper solution
for $W_k$, while $Z$ is a lower solution
for $v_k$.

Since $u_0=0$ on $D$, we have
$W_k(0)=0\le\Phi$ and $Z(0)=b_0\le v_0$.
For each fixed $k>0$, the bound
$W_k\le k^{1/m_3}\|u_0\|_\infty$
and the divergence of $\Phi$ as $|x|\uparrow R$
give $W_k\le\Phi$ on
$\partial B_r\times[0,T]$
for $r<R$ sufficiently close to $R$.
On the same boundary, \eqref{eq:vformula} gives
\[
\begin{aligned}
v_k(x,t)^{-q}
&=v_0(x)^{-q}
+q\int_0^t W_k(x,\tau)^{m_3}\,d\tau\\
&\le b_0^{-q}+qt\Phi(x)^{m_3}
=Z(x,t)^{-q}.
\end{aligned}
\]
Hence $v_k\ge Z$ on this boundary.
Applying the comparison principle on $B_r$
and letting $r\uparrow R$, we obtain
\begin{equation}\label{eq:starbounds}
W_k\le\Phi,\qquad v_k\ge Z>0
\quad\text{in }B_R\times[0,T].
\end{equation}
Thus, for every $K\Subset D$ and $T>0$,
we have $W_k\le C$ and $v_k\ge c>0$
on $K\times[0,T]$,
where $C,c$ are independent of $k$.

By monotonicity, the limits
\[
W_*=\sup_{k>0}W_k,\qquad
V_*=\inf_{k>0}v_k
\]
exist, with $W_*<\infty$ and $V_*>0$
on $D\times[0,T]$.
Since $W_k^{m_3+1}v_k^{m_2}$ is uniformly bounded
on compact subsets of $D\times(0,\infty)$,
interior parabolic estimates, the parabolic
Sobolev embedding, and monotonicity give
$W_k\to W_*$ locally uniformly on
$D\times(0,\infty)$.

To extend the convergence to $t=0$,
fix $K\Subset G\Subset D$, where $G$ has
smooth boundary.
By \eqref{eq:starbounds},
$W_k\le C$ on $\overline G\times[0,T]$.
Since $W_k(0)=0$ in $G$ and
$\partial_t W_k\le\Delta W_k$,
comparison with the solution of the heat equation
in $G$ with zero initial value and boundary value $C$
gives
\[
\lim_{t\downarrow0}\sup_{k>0}
\|W_k(t)\|_{L^\infty(K)}=0.
\]
Thus $W_k\to W_*$ locally uniformly on
$D\times[0,T]$.

On $D$, the second equation of \eqref{eq:scaled}
gives
\[
v_k(x,t)=
\left[v_0(x)^{-q}
+q\int_0^t W_k(x,s)^{m_3}\,ds\right]^{-1/q}.
\]
The local uniform convergence of $W_k$
gives \eqref{eq:Vstar} and
$v_k\to V_*$ locally uniformly on
$D\times[0,T]$.
Passing to the limit in \eqref{eq:scaled},
we obtain \eqref{eq:starsystem}, with
$W_*(0)=0$ and $V_*(0)=v_0$ in $D$.

The comparison argument leading to
\eqref{eq:scalarlower} also applies when
$\delta=1$.
In this case, $W_k\ge w$ for $k\ge1$,
where $w$ solves \eqref{eq:scalar}
and is positive in $\Omega\times(0,\infty)$.
Hence $W_*>0$ on $D\times(0,\infty)$, and
\[
0<V_*(x,t)=
\left[v_0(x)^{-q}
+q\int_0^t W_*(x,s)^{m_3}\,ds\right]^{-1/q}
<v_0(x)
\]
for $x\in D$ and $t>0$.

Finally, extend $V_*$ by zero to $\Omega\setminus D$.
Since $0\le v_k,V_*\le v_0$ and
$v_k=V_*=0$ outside $D$,
the local uniform convergence and exhaustion
by compact subsets of $D$ give
$v_k\to V_*$ in $C([0,T];L^r(\Omega))$
for every finite $r$. This
completes the proof.
\end{proof}

\section{The critical regime}\label{sec:critical}

Suppose $\sigma=0$, so $m_4=m_2+1$.
By \eqref{eq:kmonotone}, $u_k$ decreases to
a limit $u$, and $A_k\downarrow A_u$.
Together with Lemma~\ref{lem:compact}, this gives
$u_k\to u$ in $C([0,T];L^r(\Omega))$
for every finite $r$.

We first consider $\delta\le0$.

\begin{proof}[Proof of Theorem~\ref{thm:classification}(iii)
for $\sigma=0$ and $\delta\le0$]
On $D$, we have $v_k>0$. Since $m_4=m_2+1$,
the second equation of \eqref{eq:four} gives
\[
\partial_t\log\frac{v_0}{v_k}
=ku_k^{m_3}v_k^{m_2}
=u_k^{-\delta}R_k
\le\|u_0\|_\infty^{-\delta}R_k.
\]
By \eqref{eq:mass},
\begin{equation}\label{eq:logbound}
\int_D\log\frac{v_0(x)}{v_k(x,T)}\,dx
\le\|u_0\|_\infty^{-\delta}\|u_0\|_1.
\end{equation}
On the other hand, \eqref{eq:vformula} gives
\[
\log\frac{v_0(x)}{v_k(x,T)}
=\frac1{m_2}\log\bigl(1+m_2kv_0(x)^{m_2}A_k(x,T)\bigr).
\]
Since $A_k\downarrow A_u$, the right-hand side
diverges wherever $x\in D$ and $A_u(x,T)>0$.
Fatou's lemma and \eqref{eq:logbound} give
$A_u(x,T)=0$ for almost every $x\in D$.
Thus $u=0$ a.e. on $D\times(0,T)$,
and Lemma~\ref{lem:fixed} gives $u=U_D$
and \eqref{eq:localzero}.

Fix $K\Subset D$ and choose
$0\le\psi\in C_c^\infty(D)$ with $\psi=1$ on $K$.
Since $u_0\psi=0$ and $u=0$ on $D\times(0,T)$,
integration of the first equation in \eqref{eq:four}
gives
\[
\begin{aligned}
\int_0^T\!\!\int_\Omega\psi R_k\,dx\,dt
&=-\int_\Omega\psi u_k(T)\,dx
+\int_0^T\!\!\int_\Omega u_k\Delta\psi\,dx\,dt
\longrightarrow0.
\end{aligned}
\]
Moreover,
\[
-\partial_tv_k=u_k^{-\delta}v_kR_k
\le\|u_0\|_\infty^{-\delta}\|v_0\|_\infty R_k.
\]
Hence $v_k(T)\to v_0$ in $L^1(K)$.
Since $0\le v_0-v_k\le v_0$ and
$v_k=v_0=0$ outside $D$,
exhaustion by compact subsets of $D$ gives
$v_k(T)\to v_0$ in $L^1(\Omega)$.
By boundedness and time monotonicity,
$v_k\to v_0$ in $C([0,T];L^r(\Omega))$
for every finite $r$. This completes the
proof.
\end{proof}

Suppose $\delta>0$.
On $D$, \eqref{eq:vformula} gives
\begin{equation}\label{eq:criticalformula}
kv_k^{m_2}=\frac1{k^{-1}v_0^{-m_2}+m_2A_k},\qquad
R_k=\frac{u_k^{m_1}}{k^{-1}v_0^{-m_2}+m_2A_k}.
\end{equation}
We first estimate $R_k$ where $A_k$ is small.

\begin{lemma}\label{lem:reaction}
Suppose $\sigma=0$ and $\delta>0$, and set
$\theta=\min\{\delta,1\}$ and $\gamma=\theta/(m_3+2)$.
For every $0<\tau<T$, there exists $C_{\tau,T}>0$
such that, for all $k>0$ and $0<h\le1$,
\begin{equation}\label{eq:smallhistory}
\int_\tau^T\!\!\int_{D\cap\{A_k\le h\}}R_k\,dx\,dt
\le C_{\tau,T}\bigl(h^{\theta/m_3}+h^\gamma\bigr).
\end{equation}
Moreover, $R_k\to\mathcal F[u]$ in $L^1(Q_T)$,
and $u$ is a weak solution of \eqref{eq:memory_problem}.
\end{lemma}

\begin{proof}
Fix $0<\tau<T$ and set
\[
J_k(x)=\int_{\tau/2}^T|\partial_tu_k(x,s)|^2\,ds.
\]
By \eqref{eq:positiveenergy},
\[
\int_\Omega J_k\,dx
=\int_{\tau/2}^T\|\partial_tu_k(t)\|_2^2\,dt
\le\frac{\|u_0\|_2^2}{2\tau}.
\]

For almost every $x\in\Omega$ and $\tau\le t\le T$,
choose $s\in[\tau/2,\tau]$ such that
$u_k(x,s)^{m_3}\le2A_k(x,t)/\tau$.
By the chain rule and the Cauchy--Schwarz inequality,
\[
\begin{aligned}
u_k(x,t)^{(m_3+2)/2}
&\le u_k(x,s)^{(m_3+2)/2}
+\frac{m_3+2}{2}
\int_s^t u_k(x,r)^{m_3/2}
|\partial_tu_k(x,r)|\,dr\\
&\le\left(\frac{2A_k(x,t)}{\tau}\right)^{(m_3+2)/(2m_3)}
+\frac{m_3+2}{2}A_k(x,t)^{1/2}J_k(x)^{1/2}.
\end{aligned}
\]
Since $0<2\gamma<1$ and $(m_3+2)\gamma=\theta$,
raising both sides to the power $2\gamma$ gives
\begin{equation}\label{eq:interpolation}
u_k^\theta\le C\bigl(
\tau^{-\theta/m_3}A_k^{\theta/m_3}
+J_k(x)^\gamma A_k^\gamma\bigr).
\end{equation}

Since $m_1=m_3+\delta$,
$\partial_tA_k=u_k^{m_3}$, and
$u_k^\delta\le C u_k^\theta$,
\eqref{eq:criticalformula} and \eqref{eq:interpolation}
give, on $D\cap\{A_k>0\}$,
\begin{equation}\label{eq:historydensity}
\begin{aligned}
R_k
&\le\frac{u_k^\delta}{m_2A_k}\,\partial_tA_k\\
&\le C\bigl(
\tau^{-\theta/m_3}A_k^{\theta/m_3-1}
+J_k(x)^\gamma A_k^{\gamma-1}\bigr)
\partial_tA_k.
\end{aligned}
\end{equation}
On $\{A_k=0\}$, both $R_k$ and $\partial_tA_k$
vanish almost everywhere.
Since $A_k$ is nondecreasing in time,
integrating \eqref{eq:historydensity} over
$\{t\in[\tau,T]:A_k(x,t)\le h\}$ gives
\[
\int_{\{t\in[\tau,T]:A_k(x,t)\le h\}}
R_k(x,t)\,dt
\le C\bigl(
\tau^{-\theta/m_3}h^{\theta/m_3}
+J_k(x)^\gamma h^\gamma\bigr).
\]
Integrating in space and using H\"older's inequality,
we obtain \eqref{eq:smallhistory}.

On $\{A_k>h\}$, \eqref{eq:criticalformula} gives
$R_k\le\|u_0\|_\infty^{m_1}/(m_2h)$.
Thus, for every measurable
$E\subset\Omega\times(\tau,T)$,
\begin{equation}\label{eq:uniformintegrability}
\int_E R_k\,dx\,dt
\le C\bigl(h^{\theta/m_3}+h^\gamma\bigr)
+\frac{\|u_0\|_\infty^{m_1}}{m_2h}|E|.
\end{equation}
Letting first $|E|\downarrow0$ and then $h\downarrow0$
gives uniform integrability of $R_k$
on $\Omega\times(\tau,T)$.

Since $u_k\downarrow u$ and $A_k\downarrow A_u$,
\eqref{eq:criticalformula} gives
\[
R_k\longrightarrow\frac{u^{m_1}}{m_2A_u}
\quad\text{a.e. on }(D \times (\tau,T))\cap\{A_u>0\}.
\]
Outside $D$, we have $R_k=0$.
Moreover, $u=0$ a.e. on $\{A_u=0\}$.
On this set, split the integral of $R_k$ into
$\{A_k\le h\}$ and $\{A_k>h\}$.
The first part is bounded by \eqref{eq:smallhistory}.
For fixed $h>0$, the second tends to zero
by dominated convergence, since
$R_k\le u_k^{m_1}/(m_2h)$ on $\{A_k>h\}$
and $u_k\to0$ a.e. on $\{A_u=0\}$.
Letting $h\downarrow0$ gives
\[
\int_{(D\times(\tau,T))\cap\{A_u=0\}}
R_k\,dx\,dt\longrightarrow0.
\]
The convergence on $\{A_u=0\}$,
\eqref{eq:uniformintegrability}, and Vitali's theorem
give
\[
R_k\to\mathcal F[u]
\quad\text{in }L^1(\Omega\times(\tau,T)).
\]

By Fatou's lemma and \eqref{eq:mass},
$\mathcal F[u]\in L^1(Q_T)$.
Combining the convergence on $\Omega\times(\tau,T)$
with \eqref{eq:trace} and the integrability of
$\mathcal F[u]$, we obtain
\[
R_k\to\mathcal F[u]\quad\text{in }L^1(Q_T).
\]

By \eqref{eq:L2}, \eqref{eq:positiveenergy},
Lemma~\ref{lem:compact}, and weak compactness,
$u$ satisfies \eqref{eq:weakclass}.
Passing to the limit in the first equation
of \eqref{eq:four} gives \eqref{eq:weakform}.
Thus $u$ is a weak solution of
\eqref{eq:memory_problem}. This completes the proof.
\end{proof}

\begin{lemma}\label{lem:maximal}
Suppose $\sigma=0$ and $\delta>0$.
Then $u$ is the greatest weak solution of
\eqref{eq:memory_problem}.
In particular, $u$ depends on $v_0$ only through $D$.
\end{lemma}

\begin{proof}
Let $w$ be any weak solution of
\eqref{eq:memory_problem}. For fixed $k>0$, set
\[
V_k[w](x,t)=
\frac{v_0(x)}{\bigl(1+m_2kv_0(x)^{m_2}
A_w(x,t)\bigr)^{1/m_2}}.
\]
Then $0\le V_k[w]\le v_0$,
$V_k[w](0)=v_0$, and
$V_k[w]\in W^{1,\infty}(0,T;L^2(\Omega))$.
Differentiation gives
\begin{equation}\label{eq:Vderivative}
\partial_tV_k[w]=-kw^{m_3}V_k[w]^{m_2+1}.
\end{equation}
Since $\partial_tA_w=w^{m_3}$,
we have $w=0$ a.e. on $\{A_w=0\}$.
On $D\cap\{A_w>0\}$,
\[
\begin{aligned}
kw^{m_1}V_k[w]^{m_2}
&=\frac{w^{m_1}}{k^{-1}v_0^{-m_2}+m_2A_w}\\
&\le\frac{w^{m_1}}{m_2A_w}=\mathcal F[w].
\end{aligned}
\]
Both sides also vanish a.e. on $\{A_w=0\}$
and outside $D$.
Hence \eqref{eq:weakform} gives
\[
\partial_tw-\Delta w+kw^{m_1}V_k[w]^{m_2}\le0
\]
in the weak sense.

Set
\[
\xi_k=(w-u_k)_+,\qquad
\eta_k=(v_k-V_k[w])_+.
\]
From the equations for $u_k$, $v_k$, and $V_k[w]$
and the inequality for $w$, we obtain
\[
\begin{aligned}
\partial_t(w-u_k)-\Delta(w-u_k)
&\le k\bigl(u_k^{m_1}v_k^{m_2}
-w^{m_1}V_k[w]^{m_2}\bigr),\\
\partial_t(v_k-V_k[w])
&=k\bigl(w^{m_3}V_k[w]^{m_2+1}
-u_k^{m_3}v_k^{m_2+1}\bigr).
\end{aligned}
\]
On $\{\xi_k>0\}$, we have $w>u_k$, while on
$\{\eta_k>0\}$, we have $v_k>V_k[w]$.
Since all four functions are bounded and the
nonlinearities are Lipschitz on bounded sets,
\[
\begin{aligned}
\bigl(u_k^{m_1}v_k^{m_2}
-w^{m_1}V_k[w]^{m_2}\bigr)\xi_k
&\le C\xi_k\eta_k,\\
\bigl(w^{m_3}V_k[w]^{m_2+1}
-u_k^{m_3}v_k^{m_2+1}\bigr)\eta_k
&\le C\xi_k\eta_k.
\end{aligned}
\]

For $0<s<t<T$, multiplying the first difference
inequality by $\xi_k$ and the second by $\eta_k$,
and integrating over $\Omega\times(s,t)$
after smoothing if necessary, we obtain
\[
\begin{aligned}
&\frac12\bigl(\|\xi_k(t)\|_2^2
+\|\eta_k(t)\|_2^2\bigr)
+\int_s^t\|\nabla\xi_k(r)\|_2^2\,dr\\
\le&
\frac12\bigl(\|\xi_k(s)\|_2^2
+\|\eta_k(s)\|_2^2\bigr)
+C_k\int_s^t
\bigl(\|\xi_k(r)\|_2^2
+\|\eta_k(r)\|_2^2\bigr)\,dr,
\end{aligned}
\]
where $C_k$ is independent of $s$ and $t$.
By Gronwall's inequality,
\[
\|\xi_k(t)\|_2^2+\|\eta_k(t)\|_2^2
\le e^{2C_k(t-s)}
\bigl(\|\xi_k(s)\|_2^2+\|\eta_k(s)\|_2^2\bigr).
\]
Since $w(0)=u_k(0)=u_0$ and
$V_k[w](0)=v_k(0)=v_0$ in $L^2(\Omega)$,
we have
$\|\xi_k(s)\|_2+\|\eta_k(s)\|_2\to0$
as $s\downarrow0$.
Thus $\xi_k=\eta_k=0$, and hence
\[
w\le u_k,\qquad v_k\le V_k[w].
\]
Letting $k\to\infty$ gives $w\le u$.
Since $u$ is a weak solution by
Lemma~\ref{lem:reaction},
it is the greatest weak solution.
Finally, \eqref{eq:memory_problem} depends on
$v_0$ only through $D$, so the same holds for $u$. This completes the proof.
\end{proof}

We next examine the penetration of $u$ into $D$.
\begin{proof}[Proof of Theorem~\ref{thm:penetration}]
The convergence of $R_k$ and the dependence of $u$
on $v_0$ only through $D$ follow from
Lemmas~\ref{lem:reaction} and~\ref{lem:maximal},
respectively.
Since $m_4=m_2+1$ and $A_k\ge A_u$,
\eqref{eq:vformula} gives
\[
0\le v_k(x,t)\le
\bigl(m_2kA_u(x,t)\bigr)^{-1/m_2}\longrightarrow0
\]
whenever $x\in D$ and $A_u(x,t)>0$.
This also proves the assertion on $v_k$ in
Theorem~\ref{thm:classification}(iv).

Suppose $0<\delta<1$.
Fix $B_R(x_0)\Subset D$ and $T>0$ such that
$C_\delta\sqrt T<R/2$,
where $C_\delta$ is given by \eqref{eq:Cdelta}.
Set
\[
b_0=\min_{\overline{B_R(x_0)}}v_0>0,
\qquad
\mu_\delta=\frac2{1-\delta}.
\]
Choose $\zeta>0$ sufficiently small that
\[
L_\zeta=\sqrt{T(1+\zeta)}\,C_\delta<R/2,
\]
and set
\[
\varepsilon_k=
\left(\frac{b_0^{-m_2}}{\zeta m_2Tk}\right)^{1/m_3},
\qquad
r_0=R-L_\zeta>R/2.
\]
Define
\[
\Phi(x)=\|u_0\|_\infty
\left(\frac{(|x-x_0|-r_0)_+}{L_\zeta}\right)^{\mu_\delta},
\qquad
S_k(x)=\Phi(x)+\varepsilon_k,
\]
and
\[
Z_k(x,t)=
\left[b_0^{-m_2}+m_2ktS_k(x)^{m_3}\right]^{-1/m_2}.
\]
Since $\mu_\delta>2$, the function $\Phi$ is $C^2$
and vanishes on $B_{r_0}(x_0)$.
For $r_0<|x-x_0|<R$, set
$h=(|x-x_0|-r_0)/L_\zeta$.
Using $L_\zeta/r_0<1$, $0<h<1$,
$\mu_\delta-2=\mu_\delta\delta$,
and \eqref{eq:Cdelta}, we obtain
\[
\begin{aligned}
\Delta\Phi
&=\frac{\|u_0\|_\infty}{L_\zeta^2}
\left[
\mu_\delta(\mu_\delta-1)h^{\mu_\delta-2}
+\mu_\delta(N-1)\frac{L_\zeta}{|x-x_0|}
h^{\mu_\delta-1}
\right]\\
&\le\frac{\|u_0\|_\infty}{L_\zeta^2}
\mu_\delta(\mu_\delta+N-2)h^{\mu_\delta-2}\\
&=\frac{\Phi^\delta}{(1+\zeta)m_2T}.
\end{aligned}
\]
The same inequality holds in $B_{r_0}(x_0)$.
Since $S_k\ge\varepsilon_k$,
\[
\frac{b_0^{-m_2}}{kS_k^{m_3}}
\le\zeta m_2T.
\]
Thus, for $0\le t\le T$,
\[
\begin{aligned}
kS_k^{m_1}Z_k^{m_2}
&=\frac{S_k^\delta}
{m_2t+b_0^{-m_2}/(kS_k^{m_3})}\\
&\ge\frac{S_k^\delta}{(1+\zeta)m_2T}
\ge\frac{\Phi^\delta}{(1+\zeta)m_2T}
\ge\Delta S_k.
\end{aligned}
\]
Moreover,
$\partial_tZ_k=-kS_k^{m_3}Z_k^{m_2+1}$.
Thus $S_k$ is an upper solution
for $u_k$, while $Z_k$ is a lower solution
for $v_k$.

Since $u_0=0$ on $D$, we have
$u_0\le S_k$ and $Z_k(0)=b_0\le v_0$
in $B_R(x_0)$.
On $\partial B_R(x_0)\times[0,T]$,
$u_k\le\|u_0\|_\infty\le S_k$.
The comparison principle
\cite[Theorem~2.2]{Tsukamoto2026PersistenceReactionTermUnbalanced}
gives
\[
u_k\le S_k,\qquad v_k\ge Z_k
\quad\text{in }B_R(x_0)\times[0,T].
\]
Since $\Phi=0$ on $B_{r_0}(x_0)$, we obtain
\begin{equation}\label{eq:smallbarrier}
u_k\le\varepsilon_k
\quad\text{on }B_{R-L_\zeta}(x_0)\times[0,T].
\end{equation}

For each $0<t<R^2/(4C_\delta^2)$,
taking $T=t$ and then letting $k\to\infty$
and $\zeta\downarrow0$ in \eqref{eq:smallbarrier}
gives \eqref{eq:penetration}.
A finite covering by interior balls gives
the initial zero interval on every $K\Subset D$.

For any $T>0$, if the integral in
\eqref{eq:nonstationary} were zero,
then $u=A_u=0$ a.e. on $D\times(0,T)$.
Thus $\mathcal F[u]=0$ on $Q_T$, and
\eqref{eq:weakform} gives $u=H$.
This contradicts the positivity of $H$ in
$D\times(0,T)$.
Hence \eqref{eq:nonstationary} holds.

Suppose $\delta\ge1$.
Since $m_2=m_4-1$ and $m_1-m_3\ge1$,
\cite[Theorem~1.1]{Tsukamoto2026PersistenceReactionTermUnbalanced}
gives Theorem~\ref{thm:penetration}(ii),
with the same limit $u$.
Moreover, $A_u(x,\tau)>0$ on $D$, so
\eqref{eq:vformula} and dominated convergence give
$\|v_k(\tau)\|_r\to0$ for every finite $r$.
By time monotonicity,
$v_k\to0$ in $C([\tau,T];L^r(\Omega))$.
This completes the proof.
\end{proof}

We now complete the proof of the classification theorem.
\begin{proof}[Proof of Theorem~\ref{thm:classification}]
Lemma~\ref{lem:compact} applies to all exponent ranges.
Section~\ref{sec:noncritical} proves (i), (ii),
and (iii) when $\sigma\ne0$.
The case $\sigma=0$, $\delta\le0$ follows from
the first proof of this section.
For $\sigma=0$, $\delta>0$,
Lemmas~\ref{lem:reaction} and~\ref{lem:maximal},
together with Theorem~\ref{thm:penetration},
prove (iv).
This completes the proof.
\end{proof}

We finally examine the persistence of a zero region
under the Dirichlet boundary condition.
\begin{proof}[Proof of Theorem~\ref{thm:permanent-zero}]
Let $\lambda_1>0$ and $\phi_1>0$ be the first
Dirichlet eigenvalue and eigenfunction of $-\Delta$
in $\Omega$.
Since $f$ is compactly supported in $\Omega$,
there exists $a>0$ such that $f\le a\phi_1$.
By comparison with the heat equation,
\begin{equation}\label{eq:permanent-decay}
0\le u(x,t)\le H(x,t)
\le a\varepsilon e^{-\lambda_1t}\phi_1(x).
\end{equation}
Put $M=a\|\phi_1\|_\infty/\lambda_1$ and, for $T>0$, set
\[
I_T(x)=\int_0^T u(x,t)\,dt,\qquad
G_T(x)=\int_0^T\mathcal F[u](x,t)\,dt.
\]
By \eqref{eq:permanent-decay},
\[
0\le I_T\le M\varepsilon.
\]
Since $u_0=0$ in $B_{3R}(x_0)$, integrating
\eqref{eq:memory_problem} in time gives
\begin{equation}\label{eq:permanent-elliptic}
\Delta I_T=u(T)+G_T\ge G_T
\quad\text{in }B_{3R}(x_0)
\end{equation}
in the weak sense.

We first estimate $G_T$ in terms of $I_T$.
For a.e. $x\in D$ and every $L,\eta>0$,
the identity $\partial_tA_u=u^{m_3}$
and H\"older's inequality give
\[
\begin{aligned}
&(A_u(x,L)+\eta)^{\delta/m_1}
-\eta^{\delta/m_1}\\
=&\frac{\delta}{m_1}
\int_0^L
\frac{u(x,t)^{m_3}}
{(A_u(x,t)+\eta)^{m_3/m_1}}\,dt\\
\le&\frac{\delta}{m_1}L^{\delta/m_1}
\left(
m_2\int_0^L\mathcal F[u](x,t)\,dt
\right)^{m_3/m_1}.
\end{aligned}
\]
Since
\[
\int_0^L u(x,t)\,dt
\le L^{1-1/m_3}A_u(x,L)^{1/m_3},
\]
letting $\eta\downarrow0$ yields
\begin{equation}\label{eq:permanent-time}
G_L(x)
\ge c_0L^{-\delta} I_L(x)^\delta,
\qquad
c_0=\frac1{m_2}
\left(\frac{m_1}{\delta}\right)^{m_1/m_3}.
\end{equation}

Fix $T>0$. For a.e. $x\in D$ with $I_T(x)>0$,
set
\[
L=\min\left\{
T,\frac1{\lambda_1}
\log\frac{2M\varepsilon}{I_T(x)}
\right\}.
\]
If $L<T$, then \eqref{eq:permanent-decay}
and the definition of $L$ give
\[
\int_L^T u(x,t)\,dt
\le a\varepsilon\phi_1(x)
\int_L^T e^{-\lambda_1t}\,dt
\le M\varepsilon e^{-\lambda_1L}
=\frac{I_T(x)}2.
\]
If $L=T$, the integral is zero.
Thus, in either case,
\begin{equation}\label{eq:permanent-half}
\int_0^L u(x,t)\,dt
\ge\frac{I_T(x)}2.
\end{equation}
Since
\[
L\le\frac1{\lambda_1}
\log\frac{2M\varepsilon}{I_T(x)},
\]
\eqref{eq:permanent-time} and
\eqref{eq:permanent-half} give
\[
G_T(x)
\ge G_L(x)
\ge c_0L^{-\delta}
\left(\frac{I_T(x)}2\right)^\delta
\ge\frac{c_0\lambda_1^\delta}{2^\delta}
\frac{I_T(x)^\delta}
{(\log(2M\varepsilon/I_T(x)))^\delta}.
\]

Choose $\kappa=(1+\delta)/2$, so that
$\delta<\kappa<1$.
Since
$(\log(2/y))^\delta\le C y^{\delta-\kappa}$
for $0<y\le1$, the preceding estimate and
$0<I_T(x)\le M\varepsilon$ give
\begin{equation}\label{eq:permanent-coercivity}
G_T\ge c\varepsilon^{\delta-\kappa}I_T^\kappa
\quad\text{a.e. in }D,
\end{equation}
where $c>0$ is independent of $T$ and
$\varepsilon$.
The inequality also holds where $I_T=0$.

Set
\[
\mu=\frac2{1-\kappa},\qquad
\ell_\varepsilon=
\left(
\frac{\mu(\mu+N-2)M^{1-\kappa}}{c}
\right)^{1/2}
\varepsilon^{(1-\delta)/2}.
\]
Choose $\varepsilon_*>0$ such that
$\ell_\varepsilon<R$ whenever
$0<\varepsilon\le\varepsilon_*$, and define
\[
\Phi_\varepsilon(x)=M\varepsilon
\left(
\frac{(|x-x_0|-3R+\ell_\varepsilon)_+}
{\ell_\varepsilon}
\right)^\mu
\qquad (x\in B_{3R}(x_0)).
\]
Since $\mu>2$, $\mu-2=\mu\kappa$, and
$\ell_\varepsilon<R$, direct calculation gives
\[
\begin{aligned}
\Delta\Phi_\varepsilon
&\le
\frac{\mu(\mu+N-2)M\varepsilon}
{\ell_\varepsilon^2}
\left(
\frac{(|x-x_0|-3R+\ell_\varepsilon)_+}
{\ell_\varepsilon}
\right)^{\mu-2}\\
&=c\varepsilon^{\delta-\kappa}
\Phi_\varepsilon^\kappa.
\end{aligned}
\]
Moreover, $\Phi_\varepsilon=M\varepsilon\ge I_T$
on $\partial B_{3R}(x_0)$.

By \eqref{eq:permanent-elliptic} and
\eqref{eq:permanent-coercivity},
\[
-\Delta I_T+c\varepsilon^{\delta-\kappa}I_T^\kappa
\le0
\quad\text{in }B_{3R}(x_0),
\]
whereas
\[
-\Delta\Phi_\varepsilon
+c\varepsilon^{\delta-\kappa}
\Phi_\varepsilon^\kappa\ge0.
\]
Multiplying the difference of these inequalities by
$\chi=(I_T-\Phi_\varepsilon)_+
\in H_0^1(B_{3R}(x_0))$
and integrating over $B_{3R}(x_0)$, we obtain
\[
\int_{B_{3R}(x_0)}|\nabla\chi|^2\,dx
+c\varepsilon^{\delta-\kappa}
\int_{B_{3R}(x_0)}
(I_T^\kappa-\Phi_\varepsilon^\kappa)
\chi\,dx\le0.
\]
Both terms are nonnegative, so
$I_T\le\Phi_\varepsilon$.
Since $\Phi_\varepsilon=0$ on
$B_{3R-\ell_\varepsilon}(x_0)
\supset B_{2R}(x_0)$,
we have $I_T=0$ there.
As $u\ge0$ and $T>0$ is arbitrary,
$u=0$ a.e. on
$B_{2R}(x_0)\times(0,\infty)$.
This completes the proof.
\end{proof}

\section*{Acknowledgments}
This work was supported in part by JSPS KAKENHI Grant Numbers JP23H00085 and JP26K17020.


\end{document}